\documentclass[letterpaper, 10 pt, conference]{ieeeconf}
\IEEEoverridecommandlockouts                                             
  \usepackage{braket,amsfonts}
\usepackage{bm}
\usepackage{array}
\usepackage{placeins,graphicx,diagbox }
\usepackage{etex}
\usepackage{amsmath,amssymb}  
\usepackage{pgfplots}

\usepackage{amsmath,algorithm} 
\usepackage{amssymb}  
\usepackage{color}
\usepackage{cite}
\newcounter{algo}

\newtheorem{thm}{Theorem}
\newtheorem{lem}{Lemma}

\newtheorem{assum}{Assumption}
\newtheorem{prp}{Proposition}

\newtheorem{rem}{Remark}

\newtheorem{exa}{Example}

\usepackage{algorithmic}

\usepackage{graphicx,epstopdf}

\usepackage{amsopn}

\usepackage{xspace}
\usepackage{bold-extra}
\usepackage{tcolorbox}
\usepackage{placeins,diagbox }

\usepackage{wrapfig} 
\usepackage{tikz,pgfplots,mathtools,float}
\usepackage{booktabs,caption}
\usepackage[flushleft]{threeparttable}
\usepackage{subfigure}
\usepackage{multirow}
\usepackage{amssymb,latexsym}  
\usepackage{extarrows}
\usepackage{dsfont} 
\usepackage{amssymb}  
\usepackage{color}
\usepackage{cite}

\colorlet{texcscolor}{blue!50!black}
\colorlet{texemcolor}{red!70!black}
\colorlet{texpreamble}{red!70!black}
\colorlet{codebackground}{black!25!white!25}

\begin{document}

\title{\LARGE \bf
Distributed Optimization Over Markovian Switching Random Network}

\author{Peng Yi\textsuperscript{a,b}, Li Li\textsuperscript{a,b}
\thanks{ \noindent\textsuperscript{a}  The Department of Control Science and Engineering,
Tongji University,  Shanghai,   201804, China;}
\thanks{\noindent\textsuperscript{b} Shanghai Institute of Intelligent Science and Technology, Tongji University, Shanghai, China}
\thanks{ Emails: yipeng@tongji.edu.cn, lili@tongji.edu.cn}
\thanks{The work  was partially supported by the Key Research and Development Project of National Ministry of Science and Technology under grant No. 2018YFB1305304 and the National Natural Science Foundation of China under Grant No. 51475334. }
}

 \date{}
\maketitle
\begin{abstract}
In this paper, we  investigate the distributed
 convex optimization problem over a multi-agent system  with Markovian switching communication networks.
The objective function is the sum of
each agent's  local objective function, which cannot be known by other agents.
The communication network is assumed to switch over a set of weight-balanced directed graphs with a Markovian property.
We propose a consensus sub-gradient algorithm with two time-scale step-sizes
to handle the uncertainty due to the Markovian switching topologies and the absence of global gradient information.
With a proper selection of step-sizes, we prove the almost sure convergence of all agents' local estimates to the  same optimal solution when the union graph of the Markovian  network' states is strongly connected and the Markovian network is irreducible. Simulations are given for illustration of the results.
\end{abstract}
\IEEEpeerreviewmaketitle

\section{Introduction}

There is an increasing  research interest in distributed optimization over multiagent systems
due to its broad applications in engineering networks, such as  distributed parameters estimation in sensor networks \cite{zhang2012distributed,lei2016primal}, resource allocation in communication networks, \cite{wei2013distributed,shi2015extra}, and optimal power flow in power grids, \cite{dall2013distributed,yi2016initialization}.
Due to the privacy of each agent's local data and the burden of data centralization, in distributed optimization problems each agent can only manipulate  its local objective function without knowing other agents' objective functions, while the global
objective function to be optimized is usually taken as  the sum of agents' local objective functions.
Many significant distributed optimization algorithms have been proposed and analyzed, including (sub)gradient
algorithms\cite{nedic2009distributed,lobel2010distributed,shi2015extra}, dual averaging algorithms\cite{duchi2011dual}, primal-dual methods\cite{lei2016primal,yi2016initialization,lei2018asymptotic},
gradient tracking methods\cite{xu2015augmented,pu2018push}.
Please refer to \cite{sayed2014adaptation, nedich2015convergence,nedic2018distributed,yang2019survey,notarstefano2019distributed} for the survey of recent developments in distributed optimization.

In distributed optimization, the agents must cooperatively find a consensual optimal solution by sharing information locally with  network neighbors, hence, communication plays a vital role in the design and analysis of  distributed optimization algorithm.
Different communication models and graph connectivity assumptions, either deterministic or stochastic, have been discussed for different algorithms
including uniformly joint strongly connected graphs \cite{nedic2009distributed,duchi2011dual},
quantized communication \cite{yi2014quantized},
random graphs \cite{xu2017convergence,lei2018asymptotic},
broadcasting \cite{nedic2010asynchronous} and gossip communication \cite{lu2011gossip}.
In fact, the practical communication networks are essentially random and stochastic due to link failure,  uncertain quantization, packet dropout or node recreation.
Random communication networks with temporal independence assumptions have been investigated in distributed optimization. \cite{nedic2010asynchronous} established the
almost sure convergence of the consensus subgradient algorithm  to an optimal point when the agents share information through independent broadcast communications.
 \cite{lobel2010distributed} provided the almost sure convergence results for distributed subgradient algorithm when the communication link failures are independent and identically distributed over time. \cite{xu2017convergence} investigated the
 asynchronous distributed gradient method with a linear convergence rate for strongly convex functions when the graph weights are independently and identically drawn from the same probability space. \cite{jakovetic2018convergence} proved the optimal  convergence rate of distributed stochastic gradient methods for  strongly convex functions over temporally independent identically distributed  random networks.  \cite{lei2018asymptotic} investigated the  asymptotic normality and efficiency of distributed primal-dual gradient algorithm for  independent and identically random communication networks.
 \cite{yi2018distributed} gave a primal-dual algorithm for distributed resource allocation, also with independent and identically random communication networks.

Nevertheless, the practical communications over multiagent systems are usually random but with temporal correlation.
Markovian switching graphs have been  adopted for
modelling the random communication with one-step temporal dependence.
For example, \cite{huang2010stochastic, matei2013convergence, li2018distributed} have investigated the performance of averaging consensus algorithm with Markovian switching communication networks,  \cite{zhang2012distributed} have considered the
distributed parameter estimation problem over Markovian switching topologies,
and \cite{xiao2009kalman} investigated the Kalman filter with Markovian packet
losses when transmitting  the measurements to the filter.
However, to the best of  our knowledge, how to achieve distributed optimization
with Markovian switching graphs is not fully investigated,
because distributed optimization is a fundamentally different task from consensus or parameter estimations, except that  \cite{lobel2011distributed} have studied distributed optimization over a switching
state-dependent graphs. We also note that \cite{alaviani2019distributed} investigated the distributed optimization through the fixed  points iteration of random operators derived from a general class of  random graphs.

Motivated by the above,  we investigate the consensus subgradient algorithm to achieve optimal consensus with Markovian switching topologies. The communication graph among the agents switches within a finite graph set following a Markovian chain. Note that  \cite{lobel2011distributed} assumed that the random link failure  is dependent on the node state rather than the previous step communication, hence, it considered a different Markovain model from the Markovain random graph considered here.
We propose to select two different step-sizes for the consensus term and the gradient term to  balance the speed of consensus and innovation.
We find a sufficient choice of step-sizes to ensure that the consensus  term is slightly ``faster" than the innovation gradient term, and then we can give a mean consensus error bounds  under the Markovian assumption.
With these error bounds, we prove that all the agents converge to the same optimal solution with probability $1$.

The paper is organized as follows.
 We give the formulation of the distributed optimization  problem
   and Markovian switching communication model in Section II.
  We give the algorithm and sketch the main results in
 Section III.  We give the proofs of main theorems with an illustrative numerical example in Section IV,
and present the conclusions in section VI.
The proof of a key lemma is given  in the Appendix.

{\it Notations}:
Denote $\mathbf{1}_m=(1,...,1)^T \in \mathbb{R}^m$ and
$\mathbf{0}_m=(0,...,0)^T \in \mathbb{R}^m$. For a column vector $x
\in \mathbb{R}^m$, $x^T$ denotes its transpose. $I_n$ denotes the identity matrix in $\mathbb{R}^{n\times
n}$.
For a matrix $A=[a_{ij}] \in \mathbb{R}^{N \times N}$, $a_{ij}$
stands for the  $(i,j)_{th}$ entry in $A$. A matrix $A$ is nonnegative if $a_{ij} \geq 0, \forall i,j =1,\cdots,N$.
A nonnegative matrix $A$ is called row stochastic iff $A\mathbf{1}_N =\mathbf{1}_N$, and column stochastic matrix iff $\mathbf{1}_N^T A =\mathbf{1}_N^T$,
while $A$ is  doubly stochastic  iff $A$ is both row and column stochastic matrix.
$\otimes$ stands for the Kronecker product of two matrixes. For a probability space $(\Xi, \mathcal{F}, \mathbb{P})$, $\Xi$ is the sample space, $\mathcal{F}$ is the $\sigma$-algebra and $\mathbb{P}$ is the probability measure.  For $k=0,1,2,\cdots,$ $(v_k,\mathcal{F}_k)$  is an adapted sequences if $\sigma(v_k)\in \mathcal{F}_k$ for all $k$.
The expectation of a
random variable is denoted as $\mathbb{E}[\cdot]$.

A directed graph  $\mathcal{G} = \{ \mathcal{V}, \mathcal{E }_{\mathcal{G }}, {A}_{\mathcal{G }}\}$ is defined with
node set $\mathcal{V}=\{1,...,N\}$,  edge set
$\mathcal{E }_{\mathcal{G }}  \subset  \mathcal{V } \times \mathcal{V } $, and adjacency matrix ${A}_{\mathcal{G }} =[a_{ij} ]\in  \mathbb{R}^{N\times N}$.
   $(j,i)\in\mathcal{E }_{\mathcal{G }}$ if and only if agent  $i$ can get information from agent $j$.
${A}_{\mathcal{G }} =[a_{ij} ]$ is nonnegative and row stochastic, and $ 0< a_{ij} \leq 1$ if $(j, i)\in \mathcal{E }_{\mathcal{G }} $, and  $a_{ij}=0$, otherwise.
  Denote by  $\mathcal{N}_i=\{j|(j,i) \in
\mathcal{E}_{\mathcal{G}}\}$ the  neighbor set of agent $i$.
A path of graph $\mathcal{G}$ is a sequence of distinct agents in
$\mathcal{V}$ such that any consecutive agents in the sequence
corresponding to an edge of the graph $\mathcal{G}$. Agent $j$ is said to be
connected to agent $i$ if there is a path from $j$ to $i$.
Graph $\mathcal{G}$ is strongly connected if any two agents are
connected.
Graph $\mathcal{G}$ is  called weighted-balanced if adjacency matrix $A$ is doubly stochastic , i.e., $\mathbf{1}^T_N A_\mathcal{G}=\mathbf{1}^T_N A^T_\mathcal{G}$.
Denote by  ${D}_{\mathcal{G }}= \textrm{diag}\{ \sum_{j=1}^N a_{1j},..., \sum_{j=1}^N a_{Nj} \}$, called the in-degree matrix  of  $\mathcal{G}$. Then, the (weighted) Laplacian  matrix of $\mathcal{G } $ is
${L }_{\mathcal{G }} :={D}_{\mathcal{G }}-{A}_{\mathcal{G }}$.
When graph $\mathcal{G}$  is   strongly connected, 0 is a simple
eigenvalue of Laplacian ${L}_{\mathcal{G}}$ with the eigenspace $\{ \alpha \mathbf{1}_N| \alpha\in \mathbb{R}\}$.

\section{Problem Formulation}

In this section, we formulate the distributed optimization problem.

Consider a multi-agent network with agent (node) set $\mathcal{V}=\{1,...,N\}$, where  agent $i$ has its own objective function $f_i(x)$ unknown to any other agents.
The task  is to find the optimal solution of the sum of all the local objective functions, that is,
 \begin{equation}\label{pro}
       \min_{x \in \mathbb{R}^n} \quad  f(x), f(x)=\sum_{i=1}^N f_i(x),
 \end{equation}
 where $f_i(\cdot):\mathbb{R}^n \rightarrow  \mathbb{R}$, as a lower semicontinuous ( possible nonsmooth) convex function, is the local  objective function of agent $i$,
 and  $f(\cdot)$ is the global objective function.
We give  the following assumption on the objective functions.

\begin{assum}\label{asum1}
\begin{enumerate}
  \item  The optimization problem in (\ref{pro}) is solvable, i.e., there exists a finite $x^*\in \mathbb{R}^n$ such that
\begin{equation}
x^* \in X^* \triangleq \arg \min f(x), f(x)=\sum_{i=1}^N f_i(x) \nonumber
\end{equation}
  \item   The sub-gradient sets of $f_i(x),$ are uniformly bounded for all $i\in \mathcal{V}$, i.e.,  there exists a constant $l$  such that
$\forall g(x)\in \partial f_i(x)$,  $\parallel g(x) \parallel \leq l$, $\forall x\in \textrm{dom}(f_i)$, $\forall i \in \mathcal {V}$.
\end{enumerate}
\end{assum}

\vskip 3mm

We assume the agents exchange information locally through
a Markovian switching random communication network.
All the possible communication topologies form a set of a finite number of graphs:  $\{\mathcal{G}_1, \cdots, \mathcal{G}_m\}$ with each graph endowed with an adjacency matrix $A_{\mathcal{G}_i}$.
The time is slotted as $k=1,2,\cdots,$.
And then, we use a  random process $\theta(k)$, which is  a Markovian chain  on a finite index set $\mathcal{I}=\{1,...,m\}$
with a stationary transition matrix $P=[p_{ij}] \in \mathbb{R}^{m \times m}$,  to indicate the communication graph  at time $k$, i.e.,  $\mathcal{G}(k)=\mathcal{G}_{i}$ when $\theta(k)=i$.
The markovian property of  $\theta(k)$ implies that given the graph at time $k$  being $\mathcal{G}_{i}$, the probability of  the communication graph at time $k+1$ being $\mathcal{G}_j$ is $p_{ij}$.
The works about average consensus in \cite{huang2010stochastic,zhang2012distributed,matei2013convergence} have provided detailed descriptions and motivations for
 using  Markovian switching communication networks in distributed computation over multi-agent systems, including wireless sensor networks and UAV swarms.

Here is the assumption on the Markovian communication graphs:
\begin{assum}\label{asum2}
\begin{enumerate}
\item The adjacency matrixes ${A}_{\mathcal{G}_i}$ of each  graph
in the set $\{\mathcal{G}_1, \cdots, \mathcal{G}_m\}$  is a doubly stochastic matrix, and the union graph
$$\mathcal{G}_c \triangleq
  \bigcup_{i=1} ^m \mathcal{G}_i =\{\mathcal{V }, \bigcup_{i=1} ^m \mathcal{E }_{\mathcal{G }_i},
   \frac{1}{m}\sum_{i=1}^m {A}_{\mathcal{G }_i}\}$$ is strongly connected.

\item  The Markovian chain $\theta(k)$ is irreducible.
\end{enumerate}
\end{assum}

\section{ Distributed Algorithm and main results}

In this section, we provide the algorithm with the main results.

Denote by $x_i(k)\in \mathbb{R}^n $ the estimate of agent $i$ for the optimal solution $x^*$ at time $k$. The random variable $\theta(k)$ evolves as a markovian chain.
The communication graph takes $\mathcal{G}(k) \triangleq  \mathcal{G}_{\theta(k)}=(\mathcal{V},\mathcal{E}_{ \mathcal{G}_{\theta(k)} },A_{\mathcal{G}_{\theta(k)}})$ at time $k$.
Agent $i$ can get  the estimates of its  neighboring agents
 $\mathcal{N}_i(k)=\{j| (j,i) \in \mathcal{E}_{\mathcal{G}_{\theta(k)}}\}$ with $\mathcal{G}_{\theta(k)}$.
And then, each agent updates its  estimate with  the following algorithm
\begin{algorithm}[htbp]
\caption{Consensus subgradient algorithm}  \label{alg_d}
 {\it Initialize:} Agent $i\in \mathcal{V}$ picks an initial state ${x}_{i}(0)\in \mathbb{R}^n$.

{\it Iterate until convergence}\\
\noindent At  time $k$,  each agent $i\in \mathcal{V}$ gets its neighbour states  $\big\{{x}_j(k) \big\}_{j\in \mathcal{N}_i(k)}$ through the random graph $\mathcal{G}(k)$,  and updates its local state as follows
\begin{equation}\label{update}
\begin{split}
x_i(k+1)=x_i(k) &+ \alpha_k \sum_{j=1} ^N a_{ij}(k)(x_j(k)-x_i(k)) \\
               & - \beta_k d_i(k),
\end{split}
\end{equation}
where $ \alpha_k>0 $ and $\beta_k>0$ are the step-sizes, $a_{ij}(k)$ is $(i,j)_{th}$ entry  of ${A}_{\mathcal{G}_{\theta(k)}}$,
and $d_i(k) \in  \partial f_i(x_i(k))$ is a (sub)gradient vector of $f_i(x)$ at $x_i(k)$.
\end{algorithm}

Algorithm \ref{alg_d} is an extension of the distributed subgradient algorithm in \cite{nedic2009distributed},\cite{lobel2010distributed} by adding an additional step-size.
In equation (\ref{update}), the first consensus term drives each agent's state towards
the averaging of all agents' states, while the second term provides
the innovative gradient information to search for the optimal solution $x^*$.

To guarantee the algorithm convergence even with a randomly switching network, we have two different step-sizes $\alpha_k$ and $\beta_k$
to control the speed of consensus and innovation.
In fact, we require that ``consensus" speed is a bit of faster than ``innovation" term as specified by the following assumption.
\begin{assum}\label{asum3}
 We take the step-sizes in \eqref{update} as
 \begin{equation}
 \alpha_k=\frac{a_1}{(k+1)^{\delta _1}}, \quad
  \beta_k=\frac{a_2}{(k+1)^{\delta_2}},
\end{equation}
 where $a_1>0$, $a_2>0$, $0< \delta_1 <\delta_2 \leq 1$,
and $ \delta_2-\delta_1 \geq \frac{1}{2}.$
\end{assum}

\vskip 3mm
Now we are ready to present the main analysis results for Algorithm \ref{alg_d}.

\begin{thm}[Almost sure consensus]\label{thm1}
 Suppose Assumptions \ref{asum1}, \ref{asum2} and \ref{asum3} hold. Let $x_{i}(k), i \in \mathcal{V}$ be generated by  (\ref{update}), and  $y(k) =\frac{1}{N}\sum_{i=1}^N x_i(k)$. Then the following statements hold.

 1) The agents' states reach consensus and track the averaging of all the agents' states asymptotically with probability $1$, i.e.,
 \begin{equation}
  \lim_{k \rightarrow \infty}  \parallel x_i(k)-y(k)\parallel=0, \quad \forall i \in \mathcal{V},\quad a.s.
  \end{equation}

 2) The accumulation of the norm of track error  $y(k)-x_i(k)$ weighted by the step-sizes $\beta_k$ is bounded for each agent, i.e.,
  \begin{equation}
 \sum_{k=1}^{\infty} \beta_k \parallel y(k)-x_i(k)\parallel < \infty,  \quad \forall i \in \mathcal{V},\quad a.s.
 \end{equation}

\end{thm}

\vskip 3mm

\begin{rem}\label{rem1}
Theorem \ref{thm1} shows that all the agents almost surely reach consensus asymptotically. In fact,  we can also
show the convergence rate for reaching consensus is dominated by the difference between $\delta_1,\delta_2$. Specifically,
, we can find a $\tau < \delta_2-\frac{1}{2}$ such that
$$  \lim_{k \rightarrow \infty} (k+1)^{\tau}  \parallel y(k) -x_i(k)\parallel=0 \quad a.s., \forall i \in \mathcal{V}$$
\end{rem}

\vskip 3mm

\begin{thm}[Almost sure converge to a consensual  solution]\label{thm2}
Suppose Assumptions \ref{asum1},\ref{asum2} and \ref{asum3} hold.
 Then with Algorithm \ref{alg_d}, all the agents' states almost surely converge to the same optimal solution of \eqref{pro}, i.e.,
 $$ \lim_{k\rightarrow \infty} x_i(k) = x^*, \forall i\in \mathcal{V},  \quad a.s.$$
\end{thm}

\section{The proof of main results}
In this section, we give the proofs of the main results.
\vskip 3mm
Denote $X(k)= (x_1^T(k), \cdots, x_N^T(k))^T \in \mathbb{R}^{nN}$ and
$d(k)=(d_1^T(k), \cdots, d_N^T(k))^T \in \mathbb{R}^{nN} $,
and we can  rewrite the overall updating  equations in Algorithm \ref{alg_d}
in a compact form as
\begin{equation}\label{vector}
   X(k+1)=X(k) - \alpha_k ({L}_{\mathcal{G}_{\theta(k)}} \otimes I_n) X(k)  -\beta_k d(k),
\end{equation}
where ${L}_{\mathcal{G}_{\theta(k)}}$  is the Laplacian of  $\mathcal{G}_{\theta(k)}$.
With abuse of notation, we also use $A(k)$ to denote the random matrix ${A}_{\mathcal{G}_{\theta(k)}}$ for simplicity. Since $A(k)$ is doubly stochastic,
 (\ref{vector}) can also be written as:
\begin{equation}\label{vector2}
 X(k+1)=X(k) +\alpha_k((A(k)-I_N) \otimes I_n)X(k) -\beta_k d(k).
 \end{equation}

\vskip 3mm
To investigate the consensus of all agents' states, we define
 $$Q= \left(
       \begin{array}{c}
         Q_1 \\
          \frac{\mathbf{1}_N^T}{\sqrt{N}} \\
       \end{array}
     \right),
$$ with $Q_1 \mathbf{1}_N=\mathbf{0}$ and  $Q_1Q_1^T = I_{N-1}$. Then we have $QQ^T=I_N
     $,  i.e. $Q$ is an orthogonal matrix.
Denote  $$\Gamma= I_N -\frac{\mathbf{1}_N\mathbf{1}_N^T}{N} $$ as the disagreement matrix.
     Since $A(k)$ is  a doubly  stochastic matrix,
     $$ Q A(k) Q^T=\left(
                     \begin{array}{cc}
                       Q_1 A(k) Q_1^T & \mathbf{0}_{N-1} \\
                      \mathbf{0}_{N-1}  & 1 \\
                     \end{array}
                   \right),
       Q\Gamma=\left(
       \begin{array}{c}
         Q_1 \\
          \mathbf{0} \\
       \end{array}
     \right),
     $$ and therefore,
     \begin{equation}
     \begin{split}
      Q\Gamma A(k)  &= \left(
       \begin{array}{c}
         Q_1 \\
         \mathbf{0} \\
       \end{array}
     \right) Q^T   \left(
                     \begin{array}{cc}
                       Q_1 A(k) Q_1^T & \mathbf{0} \\
                      \mathbf{0}  & 1 \\
                     \end{array}
                   \right)  \left(
       \begin{array}{c}
         Q_1 \\
          \frac{\mathbf{1}_N^T}{\sqrt{N}} \\
       \end{array}
     \right) \\
     &=  \left(
                     \begin{array}{cc}
                       Q_1 A(k) Q_1^T Q_1  \\
                      \mathbf{0}   \\
                     \end{array}
                   \right) .
     \end{split}\nonumber
     \end{equation}

      Therefore, by multiplying both sides of  (\ref{vector2})  with $ Q\Gamma \otimes I_n$,
     \begin{equation}\label{vector3}
     \begin{split}
   & \left( \begin{array}{c}
         Q_1\\
          \mathbf{0} \\
       \end{array}
     \right) \otimes I_n X(k+1)=    \left( \begin{array}{c}
         Q_1 \\
          \mathbf{0} \\
       \end{array}
     \right) \otimes I_n  X(k)
     \\ & \qquad+ \alpha_k  \left(
                     \begin{array}{cc}
                      ( Q_1 A(k) Q_1^T -I_{N-1}) Q_1   \\
                      \mathbf{0}   \\
                     \end{array}
                   \right)\otimes I_n
                   X(k) \\
     & \qquad - \beta_k \left(
      \begin{array}{c}
          Q_1 \\
               \mathbf{ 0}\\
                      \end{array}
              \right)\otimes I_n
                  d(k).
 \end{split}
\end{equation}
     Denote $\xi(k)=   (Q_1 \otimes I_n )  X(k)$ and $H(k)=Q_1 A(k) Q_1^T -I_{N-1}$, and we have the reduced recursion of \eqref{vector3} as
     \begin{equation}\label{error}
      \xi(k+1)=\xi(k)+\alpha_k (H(k)\otimes I_n ) \xi(k) -\beta_k(   Q_1 \otimes I_n )d(k).
     \end{equation}

     Define a state transfer matrix
     $\Phi(k,s)$ for $k\geq s $  as
     \begin{equation}\label{state_transfer}
     \Phi(k,s)=(I_{N-1} +\alpha_k  H(k)) \cdots (I_{N-1} +\alpha_{s}  H(s)),
     \end{equation}
     and
     $\Phi(k,k+1)=I_{N-1}$.

     The state transfer matrix $\Phi(k,s)$ is a random matrix that plays a  key role in
     the convergence analysis, and we  have the following Lemma \ref{lem1}.

     \begin{lem}\label{lem1}
     Suppose Assumptions \ref{asum1}, \ref{asum2} and \ref{asum3} hold. For  the state transfer matrix $\Phi(k,s)$ defined in \eqref{state_transfer},  we have:

       (i)  There exist positive constants $c_0,c_1$ such that $$ \mathbb{E} [\parallel \Phi(k,s )\parallel ]\leq  c_0 \exp[ -c_1 \sum_{i=s} ^{k+1} \alpha_{i}], \forall k \geq s.$$

      (ii)  $\sum_{s=0} ^k \beta_{s} \mathbb{E} [\parallel \Phi(k,s+1 )\parallel] \rightarrow 0$ as $k\rightarrow \infty$.

      (iii)  $\sum_{k=0}^{\infty} \beta_{k+1} \mathbb{E}[\parallel \Phi(k,0) \parallel ]< \infty$.

      (iv) $\sum_{k=0}^{\infty} \beta_{k+1} \sum_{s=0} ^k \beta_{s} \mathbb{E}[ \parallel \Phi(k,s+1 )\parallel]< \infty$.

     \end{lem}

The proof of Lemma \ref{lem1} is given in the Appendix, which could have an independent interest since it is not related to the gradient part of the algorithm.
Next we give a martingale convergence result  from \cite{chen2006stochastic}, and a lemma for the convergence analysis.

\begin{prp}\label{prp}
Let  $(v_k,\mathcal{F}_k)$, $(\alpha_k,\mathcal{F}_k)$ be  two nonnegative adapted sequences.

(i) If $\mathbb{E}[v_{k+1}| \mathcal{F}_k] \leq v_k+\alpha_k$ and
$\mathbb{E} [\sum_{i=1}^{\infty} \alpha_i ]< \infty$,
then $v_k$ converges a.s. to a finite limit.

(ii) If $\mathbb{E}[v_{k+1}| \mathcal{F}_k] \leq v_k -\alpha_k$, then $ \sum_{i=1}^{\infty} \alpha_i < \infty, \quad a.s..$
\end{prp}

\begin{lem}\label{marg}
Let  $(v_k,\mathcal{F}_k)$, $(d_k,\mathcal{F}_k)$, and $(\alpha_k,\mathcal{F}_k)$ be  three nonnegative adapted sequences.
If $\mathbb{E}[v_{k+1}| \mathcal{F}_k] \leq v_k+\alpha_k -d_k$ and $\mathbb{E} [\sum_{i=1}^{\infty} \alpha_i ]< \infty$,
then $ \sum_{i=1}^{\infty} d_i < \infty \quad a.s.$, and $v_k$ converges a.s. to a finite limit.
\end{lem}
\textbf{Proof:} Since $\mathbb{E}[v_{k+1}| \mathcal{F}_k] \leq v_k+\alpha_k  $, and $\mathbb{E}[\sum_{i=1}^{\infty} \alpha_i] < \infty$,
 from (i) of Proposition \ref{prp}, we know that $v_k$ converges a.s. to a finite limit.

Set $u_{k+1}=v_{k+1}+\mathbb{ E}[ \sum_{i=k+1}^{\infty} \alpha_i| \mathcal{F}_{k+1}]$.
Then
 \begin{equation*}
    \begin{split}
    & \mathbb{E} [ u_{k+1} | \mathcal{F}_k] \leq   v_k+\alpha_k -d_k \\&
    \qquad + \mathbb{E}( \sum_{i=k+1}^{\infty} \alpha_i| \mathcal{F}_{k }) =u_k-d_k,
    \end{split}
 \end{equation*}
 and hence,  we have $ \sum_{i=1}^{\infty} d_i < \infty \quad a.s.$ from (ii) of Proposition \ref{prp}.
 \hfill $\Box$
\vskip 3mm
{\bf Proof of Theorem \ref{thm1}}:

1) It follows from  \eqref{error} that
\begin{equation}\label{ite}
\begin{split}
\xi(k+1)&=( \Phi(k,0) \otimes I_n) \xi(0)\\
&          -\sum_{s=0}^k \beta_s (\Phi(k,s+1)\otimes I_n)(   Q_1 \otimes I_n )d(s).
\end{split}
\end{equation}
Thus,
\begin{align}\label{ite}
  &\; \sum_{k=0}^{\infty} \beta_{k+1} \mathbb{E}[  \parallel  \xi(k+1) \parallel ]\leq \sum_{k=0}^{\infty} \beta_{k+1} \mathbb{E}[\parallel \Phi(k,0) \parallel ]  \parallel \xi(0) \parallel \nonumber \\
  &\; \qquad + \sum_{k=0}^{\infty}  \beta_{k+1} \sum_{s=0}^k \beta_s \mathbb{E} [\parallel \Phi(k,s+1)\parallel] \parallel  Q_1 \parallel \parallel d(s)\parallel \nonumber \\
  &\; \leq   \parallel \xi(0) \parallel \sum_{k=0}^{\infty} \beta_{k+1} \mathbb{E}[\parallel \Phi(k,0) \parallel] \nonumber \\
  &\; \qquad + \sqrt{N}l \parallel  Q_1 \parallel  \sum_{k=0}^{\infty} \beta_{k+1} \sum_{s=0} ^k \beta_s \mathbb{E} [\parallel \Phi(k,s+1 )\parallel]
\end{align}
Hence, with the results in Lemma \ref{lem1} and $ \parallel Q_1 \parallel = 1$,
$$  \sum_{k=0}^{\infty} \beta_{k+1}\mathbb{ E}[\parallel  \xi(k+1) \parallel ]< \infty.$$
Therefore, by the monotone convergence theorem (\cite{ash2014real}),
 \begin{equation}\label{sum1}
     \mathbb{E } [\sum_{k=0}^{\infty}  \beta_{k+1 } \parallel  \xi(k+1) \parallel ]< \infty.
 \end{equation}

 Since $ (Q\Gamma \otimes I_n) X(k) =(\xi(k)^T,    \mathbf{0})^T$  and  $Q$ is an orthogonal matrix,
 we get $\parallel (\Gamma \otimes I_n)X(k)  \parallel=  \parallel \xi(k)  \parallel $.
With \eqref{sum1},
\begin{equation}\label{res1}
   \sum_{k=1}^{\infty}  \beta_k \parallel  (\Gamma \otimes I_n)  X(k) \parallel < \infty  \quad a.s.
\end{equation}
Note that $$(\Gamma \otimes I_n) X(k) =\left(
                                   \begin{array}{c}
                                     y(k)-x_1(k) \\
                                     \vdots \\
                                    y(k)-x_N(k) \\
                                   \end{array}
                                 \right),
$$  and with \eqref{res1},  we have
\begin{equation}\label{res01}
    \sum_{k=1}^{\infty} \beta_k \parallel y(k)-x_i(k)\parallel < \infty  \quad \forall i \in \mathcal{V},  \quad a.s.
\end{equation}

2) By \eqref{update} and $y(k)=\frac{1}{N}\sum_{i=1}^n x_i(k)$,
\begin{equation}\label{ave}
  y(k+1)=y(k) -\beta_k\frac{1}{N} \sum_{j=1}^N d_j(k).
\end{equation}
Therefore, for any  $ i \in \mathcal{V}$,
\begin{equation}\label{ave2}
\begin{split}
 & y(k+1)-x_i(k+1)\\
 &=y(k) -x_i(k)+\alpha_k \sum_{j=1} ^N a_{ij}(k)(x_j(k)-x_i(k)))\\
  & \qquad - \beta_k(\frac{1}{N} \sum_{j=1}^N d_j(k)-d_i(k))\\
 &= y(k) -  \sum_{j=1} ^N \tilde{a}_{ij}(k) x_j(k)- \beta_k(\frac{1}{N} \sum_{j=1}^N d_j(k)-d_i(k)),\nonumber
  \end{split}
\end{equation}
where $\tilde{a}_{ij}(k)=\alpha_k a_{ij}(k), i \neq j$ and $\tilde{a}_{ii}(k)=1-\alpha_k$.
Note that we have $\sum_{j=1}^N \tilde{a}_{ij}(k)=1$, and therefore, with the Jensen's inequality,
\begin{equation}\label{ave3}
 \parallel y(k+1)-x_i(k+1) \parallel  \leq \sum_{j=1} ^N \tilde{a}_{ij}(k)\parallel  y(k) - x_j(k)\parallel +2 \beta_k l .\nonumber
\end{equation}

Denote by  $e_i(k)= \parallel y(k)-x_i(k) \parallel$, and $e(k)=\sum_{i=1} ^N  e_i(k)$.
Then taking square of above equations and by the convexity of  $\parallel \cdot \parallel ^2$,
\begin{equation}\label{ave4}
 e_i(k+1)^2 \leq \sum_{j=1} ^N \tilde{a}_{ij}(k)  e_j(k)^2  + 4l^2 \beta_k^2
 + 4l \beta_k\sum_{j=1} ^N \tilde{a}_{ij}(k) e_j(k). \nonumber
\end{equation}
Hence,
\begin{equation}
\begin{split}
 \sum_{i=1}^N e_i(k+1) ^2 \leq  \sum_{j=1} ^N e_j(k)^2  + 4N l^2 \beta_k^2
 + 4l \beta_k\sum_{i=1} ^N e_i(k),  \nonumber
  \end{split}
\end{equation}
and
\begin{equation}\label{ave5}
\begin{split}
\sum_{k=1}^{\infty} \sum_{i=1}^N e^2_i(k) \leq  &\sum_{i=1} ^N e_i(0)^2+ 4N l^2 \sum_{k=1}^{\infty} \beta_k^2\\
&+ 4l\sum_{k=1}^{\infty} \beta_k\sum_{i=1} ^N e_i(k).
\end{split}
\end{equation}

Because $\sum_{k=1}^{\infty } \beta_k ^2 < \infty$, $ \sum_{i=1} ^N e_i(k) ^2 $ converges with probability 1 as $k \rightarrow  \infty$ following from \eqref{res01} and \eqref{ave5}.

Meanwhile, with  \eqref{res01} and $\sum_{k=1}^{\infty}\beta(k)=\infty$, we have
$\liminf_{k \rightarrow \infty} e(k) =0$. Since   $e(k)$ is a.s. bounded,
there exists a subsequence $n_k$, such that $e(n_k) \xlongrightarrow[ k \rightarrow \infty] {}0 .$
Take any $\epsilon >0$, and then there exists a  $r$ such that
$e(n_r)<\epsilon$, $\sum_{k=n_r}^{\infty}  \beta_k e(k) <\epsilon$, $ \sum_{k=n_r}^{\infty} \beta_k^2 < \epsilon$.
Therefore,
\begin{equation} \label{con}
\begin{split}
  \sum_{i=1}^N  e_i(n_{m})^2 & \leq \sum_{i=1} ^N e_i(n_r)^2  + 4N l^2 \sum_{k=n_r}^{\infty} \beta_k^2 + 4l  \sum_{k=n_r}^{\infty}  \beta_k e(k) \\
  &\leq N\epsilon^2+ 4l(N^2+1)\epsilon \quad n_{m}=n_l+1,n_l+2 \cdots
  \end{split}
\end{equation}
Therefore, we  conclude
\begin{equation}\label{as}
\lim_{k \rightarrow \infty} e_i(k)=0, \quad a.s. \quad \forall i\in\mathcal{V}.
\end{equation}
\hfill $\Box$

Next,  we prove that all the agents almost surely converge to the same optimal solution.

{\bf Proof of Theorem \ref{thm2}:}
From \eqref{ave} it follows that, for any $ x \in X^*$,
\begin{align}\label{tm12}
&\parallel y(k+1) -x \parallel^2 = \parallel y(k) -x-\beta_k\frac{1}{N} \sum_{j=1}^N d_j(k)\parallel^2 \nonumber \\
& \qquad = \parallel y(k) -x \parallel^2 -\frac{2\beta_k}{N} \sum_{j=1}^N d_j(k)^T ( y(k) -x) \\\nonumber
& \qquad +\frac{\beta(k)^2  }{N^2} \parallel \sum_{j=1}^N d_j(k)\parallel^2. \nonumber
\end{align}

Note that
\begin{equation}
\begin{split}
  & d_j(k)^T ( y(k) -x)  \geq f_j(x_j(k)) -f_j(x) +d_j(k)^T ( y(k) -x_j(k)) \\
  &\qquad  =   f_j(x_j(k))-f_j(y(k))+f_j(y(k))-f_j(x)\\
  &\qquad +d_j(k)^T ( y(k) -x_j(k)) \\
  &\qquad \geq  -2l \parallel  y(k) -x_j(k)\parallel  +   f_j(y(k))-f_j(x). \nonumber
  \end{split}
\end{equation}
From \eqref{tm12}, we have
 \begin{equation}
\begin{split}
  &  \parallel y(k+1) -x \parallel^2 \leq  \parallel y(k) -x \parallel^2-\beta_k\frac{2}{N} [f(y(k))-f(x)] \\
  &+\beta_k \frac{4l}{N} \sum_{j=1}^N \parallel  y(k) -x_j(k)\parallel +\beta(k)^2 l^2.  \nonumber
  \end{split}
\end{equation}


By the result 1) in Theorem \ref{thm1},
$\mathbb{E} [\sum_{k=1}^{\infty} \beta_k \parallel y(k)-x_i(k)\parallel] < \infty$.
 Since $\sum_{k=0}^{\infty}\beta_k^2 < \infty$ and  $f(y(k)) \geq f^*$, the condition in
Lemma \ref{marg} is satisfied. Thus, $||y(k)-x||$ converges a.s., and
\begin{equation}\label{res3}
  \sum_{k=0}^{\infty} \beta_k  [f(y(k))-f(x)] < \infty  \quad a.s.
\end{equation}
It follows from  $\sum_{k=0}^{\infty}\beta(k)=\infty$ and \eqref{res3}  that
$$\liminf_{k \rightarrow \infty} f(y(k))-f(x)=0. \quad a.s.$$
With a similar argument of \eqref{con},  we obtain
$$\lim_{k \rightarrow \infty}  y(k) =x^*,x^* \in X^*, \quad  a.s.$$
Therefore, all the agents almost surely converge to the same optimal solution of problem (\ref{pro}).
\hfill $\Box$

\subsection{Simulation}
\begin{exa}
We give an example to illustrate the algorithm.  Consider five agents with the local objective functions as  follows:
\begin{equation}\label{sf}
\begin{array}{lll}
&\;& f_1(x)=\ln(e^{0.1x_1}+e^{0.2x_2})+5\min_{z\in \Omega}||x-z  ||; \\
&\;& f_2(x)=3(x_1)^2\ln((x_1)^2+1)+2(x_2)^2;\\
&\;& f_3(x)=3(x_1-10)^2+0.2(x_2-8)^2+ 2|x_1|+2|x_2|;\\
&\;& f_4(x)=\frac{4(x_1)^2}{\sqrt{2(x_1)^2+1}}+0.1(x_1+x_2)^2;\\
&\;& f_5(x)=(x_1+5x_2-10)^2; \\
&\;&    \qquad \quad + 4\max\{ x_1+x_2, (x_1+x_2)^2\},  \nonumber
\end{array}
\end{equation}
with a decision variable as $x=(x_1,x_2)\in \mathbb{R}^2$ and a set $\Omega$  as $\Omega=\{x \in \mathbb{R}^2 | x_1^2+x_2^2\leq 1\}$.

The five agents share information with three graphs $\{\mathcal{G}_1,\mathcal{G}_2,\mathcal{G}_3\}$,
whose weighted adjacency matrices are $A_1,A_2,A_3 \in \mathbb{R}^{5\times 5}$, respectively.
The transition matrix of the stationary Markovian chain $\theta(k)$ is $P\in \mathbb{R}^{3\times 3}$.
We let $P,A_1,A_2,A_3$ to be the following matrices:
$$\left(
                                                              \begin{array}{ccc}
                                                                0.5 & 0.5 & 0 \\
                                                                0 & 0.6 & 0.4 \\
                                                                0.2 & 0 & 0.8 \\
                                                              \end{array}
                                                            \right),
\left(
  \begin{array}{ccccc}
    0 & 0     & 0    & 0.5  & 0.5 \\
    1 & 0     & 0    & 0    & 0 \\
    0 & 0.5   & 0.5  & 0    & 0 \\
    0 & 0.5   & 0    & 0.5  & 0 \\
    0 &  0    & 0.5  & 0    & 0.5 \\
  \end{array}
\right),
$$
$$\left(
       \begin{array}{ccccc}
         0 & 0 & 0   & 0.5 & 0.5 \\
         0 & 1 & 0   & 0   & 0 \\
         1 & 0 & 0   & 0   & 0 \\
         0 & 0 & 0.5 & 0.5 & 0 \\
         0 & 0 & 0.5 & 0   & 0.5 \\
       \end{array}
     \right),
     \left(
       \begin{array}{ccccc}
         0.5   & 0 & 0     &     0  & 0.5 \\
         0     & 0.5 & 0   & 0     & 0.5    \\
         0.5     & 0.5 & 0   & 0     & 0 \\
         0     & 0 &1 & 0   & 0 \\
         0     & 0 & 0 & 1     & 0 \\
       \end{array}
     \right).
$$

We choose the step-sizes as $\alpha_k=\frac{1}{(k+1)^{0.3}}$ and $\beta_k=\frac{1}{(k+1)^{0.9}}$.
The (sub)gradients are normalized to 1.
We perform the simulation for 100 times, and the  estimations of the five agents always reach the same optimal solution. (The optimal solution
is unique in this case.)
The simulation results are shown in Figure \ref{f1} and \ref{f2}.

\begin{figure}
  \includegraphics[width=3.3in]{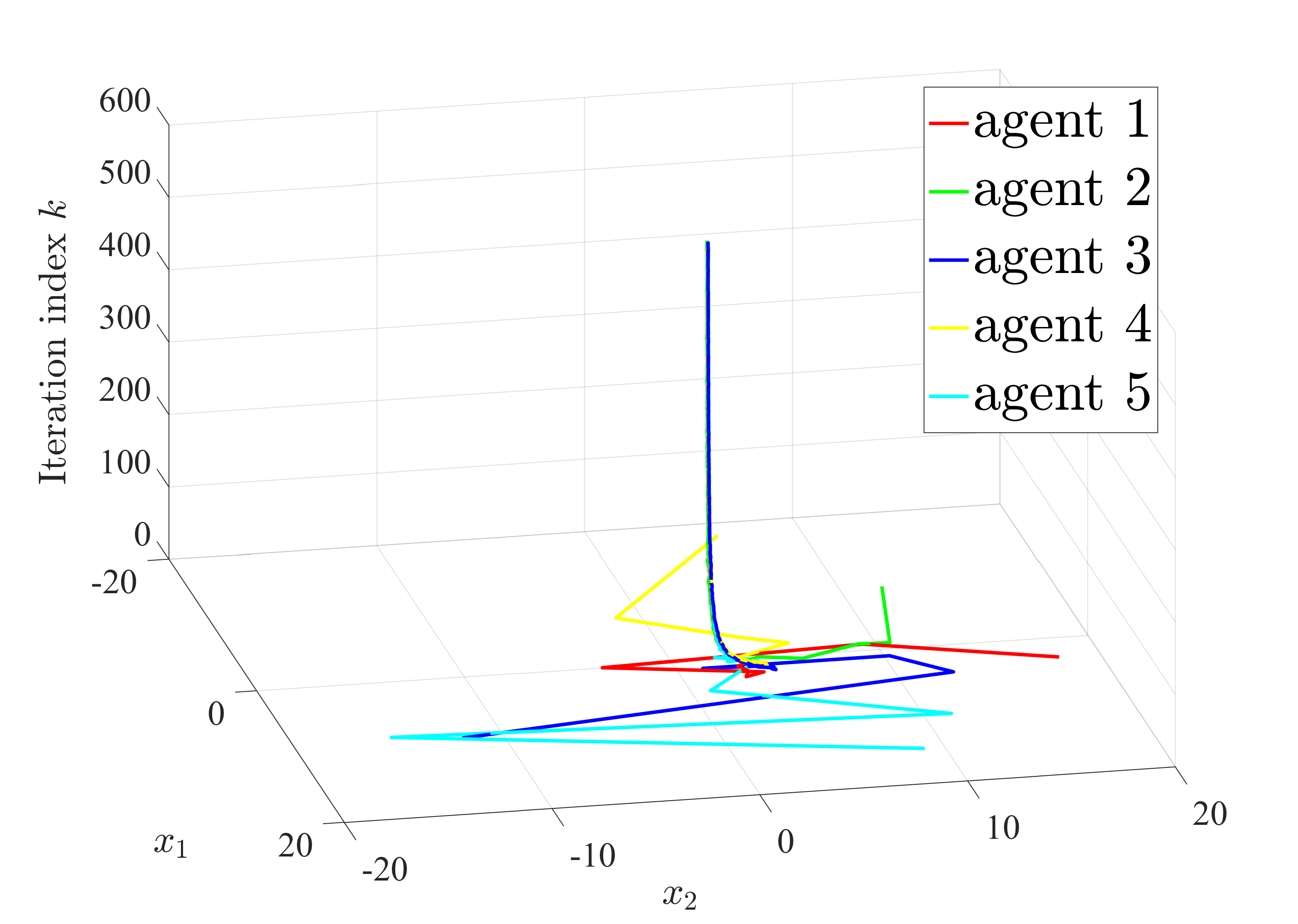}\\
  \caption{The trajectories of the five agents' states}\label{f1}
\end{figure}

\begin{figure}
  \includegraphics[width=3.3in]{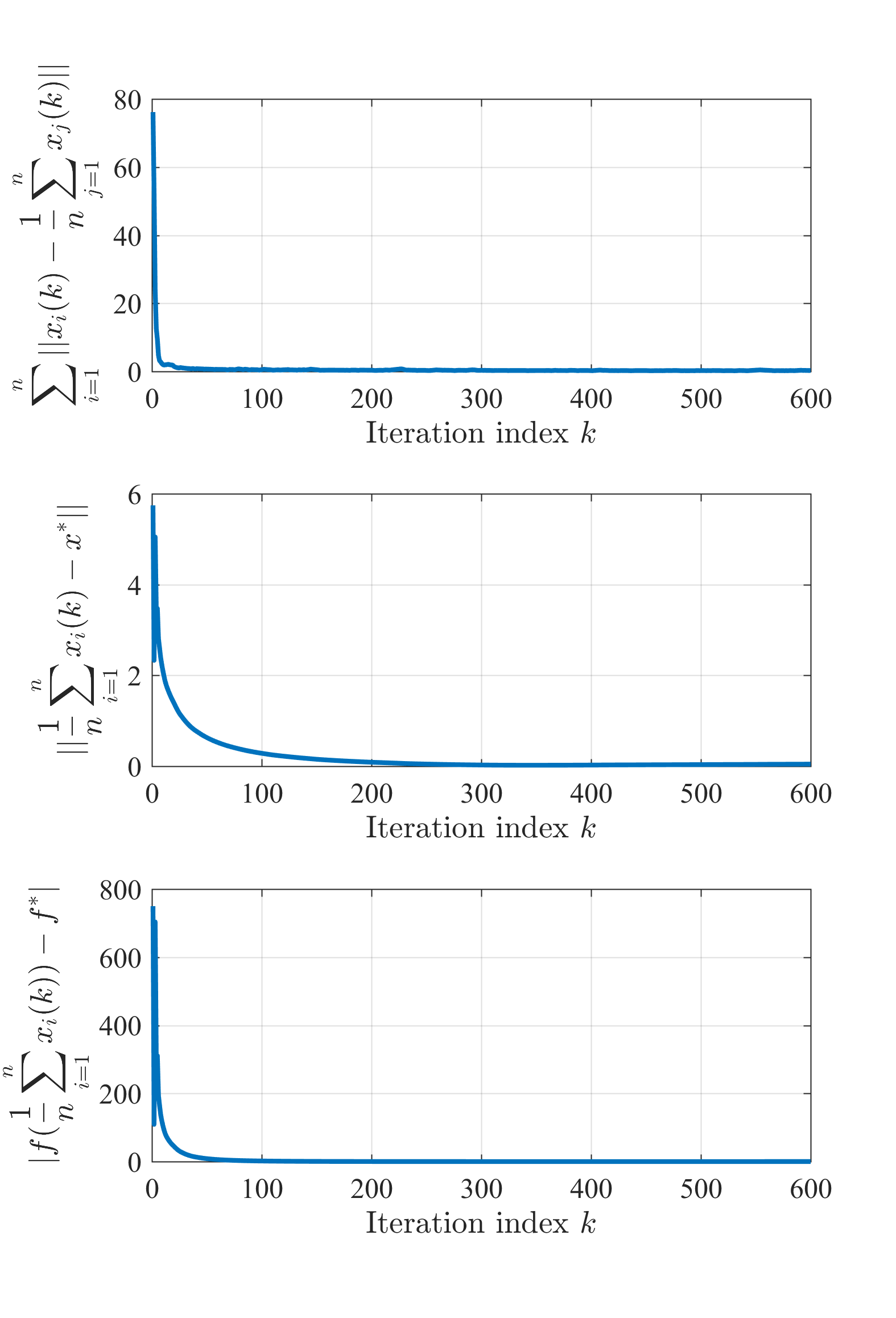}\\
  \caption{ The figure shows the trajectories of three performance index: the consensus error,  $\sum_{i=1}^5 ||x_i(k)- \frac{1}{5}\sum_{j=1}^5x_j(k) ||$, the distance to optimal solution, $|| \frac{1}{5}\sum_{i=1}^5 x_i(k)-x^* ||$ and the optimal value gap, $ |f( \frac{1}{5}\sum_{i=1}^5 x_i(k)  )-f^*  | $.}\label{f2}
\end{figure}

\end{exa}

\section{Conclusions}

In this paper, we proposed a consensus subgradient algorithm to solve a distributed optimization problem
with Markovian switching random communication networks.
The algorithm was given with two time-scale step-sizes, different from most existing ones.
We showed the almost sure convergence with a proper connectivity assumption and  step-size choices.
In the future, we will work on the mean-square convergence rate analysis.

\section*{Appendix: Proof of key lemmas}

{\bf Proof of Lemma \ref{lem1}}

(i):
Denote $\mathcal{F}_k=\sigma \{ \theta(t), 0 \leq t \leq k\} $.

     Take $h=(m-1)^2+1$ where $m$ is number of states in $\mathcal{I}$. We first prove that $\forall t$,
     $\theta(k)$ will visit all the states in $\mathcal{I}$  with a positive probability  during $[t, t+h-1]$.

     With the transition matrix $P \in \mathbb{R}^{m \times m}$ being a stochastic matrix, there exists a $\zeta >0$ such that $p_{ij} \geq \zeta $ when $p_{ij} >0$.
      Assume $\theta(t)=i_1$, and then $\forall i \in \mathcal{I}\setminus i_1$ we conclude that
      $$\mathbb{P}(\theta(k)  \; {\rm visits} \;   i \;  {\rm during}  \;   [t, t + m -1]) \geq \zeta^{m-1}$$
       because $\theta(k)$ is irreducible.
      Hence,
      $$\mathbb{P}( \theta(k) \; {\rm visits \; all  \;  states \; in \; } \mathcal{I} \; {\rm during} \;  [t, t+h-1] ) \geq  \zeta^h.$$

     Secondly, with the union graph of $\{\mathcal{G}_1,...,\mathcal{G}_m\}$ being strongly connected,
     we prove that there exists a constant $\gamma_0>0$ such that
     \begin{equation}\label{pes}
        \mathbb{E}[(\sum_{k=t} ^{t+h-1} H(k)^T+H(k)| \mathcal{F}_{t})] \leq -\gamma_0 I_{N-1}, \forall t.
     \end{equation}

     Since  $A(k)$ is a doubly stochastic matrix, we get
     $$2I_N -A(k)-A(k)^T={L}_{\mathcal{G}_{\theta(k)}}+ {L}_{\mathcal{G}_{\theta(k)}}^T=2{L}_{\mathcal{\hat{G}}_{\theta(k)}},$$
     where $\mathcal{\hat{G}}_{\theta(k)}$ is the undirected mirror graph of $\mathcal{G}_{\theta(k)}$.
     As a Laplician matrix of an undirected graph, ${L}_{\mathcal{\hat{G}}_{\theta(k)}}$  is positive semi-definite,
      i.e.,  $x^T{L}_{\mathcal{\hat{G}}_{\theta(k)}} x \geq 0, \forall x \in \mathbb{R}^N.$

     Taking  $x=Q_1^T u, u \neq \mathbf{0}$, we have
  \begin{equation}
  \begin{split}
  &u^T Q_1  (2I_N -A(k)-A(k)^T) Q_1^T u \\
  &= -u^T(H(k)+H(k)^T)u \geq 0,
  \end{split}
  \end{equation}
   and thereby,   $H(k)+H(k)^T$ is negative semi-definite.

    Similarly,
    \begin{equation}
    \begin{split}
& \mathbb{E}[u^T Q_1  (2h I_N -\sum_{k=t}^{t+h-1}[A(k)+ A(k)^T]) Q_1^T u|\mathcal{F}_t] \\
&=-\mathbb{E}[ u^T(\sum_{k=t}^{t+h-1}[H(k)+H(k)^T])u |\mathcal{F}_t]  \geq 0.
\end{split} \nonumber
\end{equation}

Therefore,
$\mathbb{E}[\sum_{k=t}^{t+h-1}[H(k)+H(k)^T]|\mathcal{F}_t]$ is also negative semi-definite.
  With Assumption \ref{asum2},  we have
     $\mathbb{E}[ x^T (2 h I_N- \sum_{k=t}^{t+h-1}[A(k) + A(k)^T])x| \mathcal{F}_t] =0$
      if and only if $  x=c \mathbf{1}_N, c \neq 0$.
      We conclude $Q_1^T u \neq c \mathbf{1}_N, \forall c \neq 0$ since otherwise  $c=0$ can also be concluded from $ \mathbf{1}_N^T  Q_1^T u = c  \mathbf{1}_N^T \mathbf{1}_N=0$.
       Therefore,
      \begin{equation}
    \mathbb{E}[u^TQ_1  (2h I_N -\sum_{k=t}^{t+h-1}[A(k)+A(k)^T]) Q_1^T u|\mathcal{F}_t] >0
     \end{equation}
     i.e., $\mathbb{E}[\sum_{k=t} ^{t+h-1} [H(k)^T+H(k)]| \mathcal{F}_{t}]$ is negative definite.
     Hence, there exists a constant $\gamma_0 >0 $ to make (\ref{pes}) hold.

As a result,
 \begin{equation}
 \begin{split}
 &\Phi(h+t-1, t) ^T \Phi(h+t-1, t)= (I_{N-1}+\alpha_{t}H(t)^T) \cdots\\
 & (I_{N-1}+\alpha_{t+h-1}H(t+h-1)^T) \\ & (I_{N-1}+\alpha_{t+h-1}H(t+h-1) )
  \cdots (I_{N-1}+\alpha_{t}H(t))\\
  & \leq I_{N-1} + \sum_{k=t} ^{t+h-1} \alpha_k (H(k)^T+H(k) )+ c_2 \alpha_t^2 I_{N-1}, \nonumber
 \end{split}
 \end{equation}
    with $c_2>0$ as a constant.
    Combined with (\ref{pes}),
\begin{equation} \label{3.3}
\begin{split}
&\mathbb{E}[ \Phi(h+t-1, t) ^T \Phi(h+t-1, t)| \mathcal{F}_{t}] \\
&\leq ( 1- \gamma_1 h \alpha_{t+h-1} + c_3 h \alpha_t^2 )I_{N-1},
\end{split}
\end{equation}
where $\gamma_1=\frac{\gamma_0 }{h }>0$, $c_3=\frac{c_2}{h }>0$.

In fact, $\forall k: t \leq k \leq t+h-1$,
\begin{equation}
\begin{split}
&( \frac{k+1}{ t+h} )^{\delta_1}=(1+ \frac{k+1-t-h}{t+h} )^{\delta_1}\\
 &=1+{\delta_1} \frac{k+1-t-h}{t+h} +  O(\frac{ h} {t+h})^2.
\end{split}
\end{equation}
Therefore, when $t$ is large enough
\begin{equation} \label{3.1}
\begin{split}
 & \alpha_k- \alpha_{t+h-1} =\frac{a_1}{(k+1)^{\delta_1}}- \frac{a_1}{(t+h )^{\delta_1}}\\
 & =\frac{a_1}{(k+1)^{\delta_1}}(1- ( \frac{k+1}{t+h} )^{\delta_1}) \\
 &  \leq \frac{a_1}{(k+1)^{\delta_1}} ( \frac{  h}{t+h} +  O(\frac{ h }{t+h})^2) \leq \frac{a_1}{4(k+1)^{\delta_1}} \\
 &= \frac{\alpha(k)}{4}
  \end{split}
\end{equation}

 We also have:
 \begin{equation}\label{3.2}
    \frac{\alpha_t^2}{ \alpha_k}=[\frac{k+1}{(t+1)^2}]^{\delta_1} \xlongrightarrow  [t \rightarrow \infty ]{ } 0.
 \end{equation}
Thereby,  with  \eqref{3.1} and \eqref{3.2}, we conclude that $\exists k_2$ when $t \geq k_2$,
 \begin{equation}
\begin{split}
&  1- \gamma_1 h \alpha_{t+h-1} + c_3 h  \alpha_t^2 \\
& \quad =1-  \gamma_1 \sum_{k=t}^{t+h-1}  \alpha_k +\gamma_1 \sum_{k=t}^{t+h-1}  (\alpha_k-   \alpha_{t+h-1})  \\
&\qquad + c_3 \sum_{k=t}^{t+h-1}  \alpha_k \frac{\alpha_t^2}{ \alpha_k}\\
& \quad \leq 1-  \gamma_1 \sum_{k=t}^{t+h-1}  \alpha_k +\frac{\gamma_1 }{4}\sum_{k=t}^{t+h-1}   \alpha_k + \frac{\gamma_1 }{4} \sum_{k=t}^{t+h-1}  \alpha_k \\
&\quad=(1-2c_4\sum_{k=t}^{t+h-1}  \alpha_k)£¬ \nonumber
  \end{split}
\end{equation}
 with a constant  $c_4= \frac{\gamma_1 }{4}>0 $.

  By \eqref{3.3},  when $t\geq k_2$,
  \begin{equation} \label{3.4}
        \mathbb{E}[ \Phi(h+t-1, t) ^T \Phi(h+t-1, t)| \mathcal{F}_{t}] \leq (1-2c_4\sum_{k=t}^{t+h-1}  \alpha_k) I_{N-1}. \nonumber
\end{equation}

Now given another integer  $ s \geq 1$, we have the following estimation by recursions,
 \begin{equation} \label{3.4}
 \begin{split}
    & \mathbb{E}[  \Phi(sh-1+t, t) ^T \Phi(sh-1+t, t) ] \\
    &=\mathbb{ E}\big[  \Phi((s-1)h-1+ t  , t )^T   \mathbb{E}[ \Phi( sh+t-1  ,(s-1)h+ t ) ^T \\
   &\Phi( sh+t-1  ,(s-1)h+ t )| \mathcal{F}_{(s-1)h+ t}]
    \Phi((s-1)h-1+ t  , t ) \big]  \\
   &\leq (1-2c_4\sum_{k=(s-1)h+ t}^{sh+ t-1}  \alpha(k))  \\
   &\mathbb{E}[  \Phi((s-1)h-1+t, t) ^T \Phi((s-1)h-1+t, t) ]\\
   & \leq \exp[-2c_4\sum_{k=   t}^{sh+ t-1} \alpha_k)]I_{N-1}  , \quad t \geq k_2,
     \end{split}
\end{equation}
based on the inequality $1-x\leq e^{ -x}, \;  \forall x \geq 0$.

Therefore,
    \begin{equation}
    \begin{split}
     \mathbb{E } [ \parallel \Phi(sh-1+t, t) \parallel] \leq c_6\exp[-c_5\sum_{k=   t}^{sh+ t-1} \alpha_k)]
    \end{split}
    \end{equation}
 with $c_6=1 $ and $c_5=2c_4$ as positive constants.

  Since  $H(k)=Q_1 {L}_{\mathcal{G}_{\theta(k)} } Q_1^T  $ and $ \mathcal{G}_{\theta(k)} $   switches
  among a finite set of graphs,  there exists a constant  $C_{max} >0$,
   such that $\parallel H(k) \parallel \leq C_{max}, \forall k \geq 1.$
Therefore,  $\forall k, s \geq k_2$, we have $\exists \varsigma \geq 0, 0 \leq  r \leq h-1$ such that $k-s=\varsigma h+r$.
Then
\begin{equation}\label{exp1}
\begin{split}
& \mathbb{E}   [ \parallel \Phi(k, s) \parallel ] \leq   \mathbb{E}  [ \parallel \Phi(k, k-\varsigma h+1) \parallel ]    \parallel \Phi(r+s, s) \parallel \\
&\leq \prod_{i=s}^{r+s}(1+\alpha_i   C_{max})  c_6  \exp[ - c_5 \sum_{i=s+r+1}^{k} \alpha_i] \\
&\leq  \prod_{i=s}^{r+s}(1+\alpha_i   C_{max})  c_6  \exp[c_5 \sum_{i=s}^{r+s}\alpha_i]\exp[ - c_5 \sum_{i=s}^{k} \alpha_i]\\
&\leq   \tilde{c}_0 \exp[ -c_5 \sum_{i=s} ^{k+1} \alpha_i],
\end{split}
\end{equation}
with $\tilde{c}_0= c_6 \prod_{i=k_2}^{k_2+h}(1+\alpha_i   C_{max}) \exp[c_5 \sum_{i=k_2}^{k_2+h}\alpha_i]$ and $c_1=c_5$ as positive constants.

When $s\leq k_2$, without loss generality we assume $k\geq k_2$, and then
\begin{equation}\label{fin}
\begin{split}
& \mathbb{E}  [\parallel \Phi(k, s) \parallel ] \leq   \mathbb{E}  [\parallel \Phi(k,k_2) \parallel ]    \parallel \Phi(k_2, s) \parallel \\
&\leq \prod_{i=s}^{k_2}(1+\alpha_i   C_{max})  \tilde{c}_0  \exp[ - c_5 \sum_{i=k_2}^{k} \alpha_i] \\
&\leq  \prod_{i=s}^{k_2}(1+\alpha_i   C_{max})  \tilde{c}_0  \exp[c_5 \sum_{i=s}^{k_2}\alpha_i]\exp[ - c_5 \sum_{i=s}^{k} \alpha_i]\\
&\leq   \hat{c}_0 \exp[ -c_1 \sum_{i=s} ^{k+1} \alpha_i],
\end{split}
\end{equation}
With $\hat{c}_0= \tilde{c}_0 \prod_{i=0}^{k_2}(1+\alpha_i   C_{max})  \exp[c_5 \sum_{i=0}^{k_2}\alpha_i]$ and $c_1=c_5$ as constants.
Taking $c_0 > \max\{ \tilde{c}_0, \hat{c}_0 \}$,  Lemma \ref{lem1} (i) holds from \eqref{exp1} and \eqref{fin}.

 (ii) Because
\begin{equation}\label{stepestimate}
\begin{split}
    &\sum_{i=s} ^{k+1} \alpha_{i}= \sum_{i=s} ^{k+1}  \frac{a_1}{(i+1)^{\delta_1}} \geq \int_{x=s}^{k+1}  \frac{a_1}{(x+1)^{\delta_1}} dx\\
    &\quad  \geq \frac{a_1}{1-\delta_1}(x+1)^{1-\delta_1}\mid {k+1 \atop s} \\
    & \quad \geq \frac{a_1}{1-\delta_1}[(k+2)^{1-\delta_1} -(s+1)^{1-\delta_1}] ,
\end{split}
\end{equation}
 we have
     \begin{equation}\label{stepcomp}
     \begin{split}
      & \sum_{s=0} ^k \beta_{s} \mathbb{E} [\parallel \Phi(k,s+1 )\parallel]  \\
      & \leq  \sum_{s=0} ^k  \frac{a_2}{(s+1)^{\delta_2}}
      c_0 \exp[\frac{a_1 c_1}{1-\delta_1}(s+2)^{1-\delta_1}\\
      & \qquad -\frac{a_1 c_1}{1-\delta_1}(k+2)^{1-\delta_1}] \\
      & \leq  a_2c_0 \frac{1}{q(k)} \sum_{s=0} ^k  \frac{1}{(s+1)^{\delta_2}}
      \exp[\frac{a_1 c_1}{1-\delta_1}(s+2)^{1-\delta_1}],
      \end{split}
     \end{equation}
with $q(k)=\exp[ \frac{a_1 c_1}{1-\delta_1}(k+2)^{1-\delta_1}]$.

     It is easy to verify that there exists a $ k_1>0$, such that $\forall x \geq k_1$,
      $\frac{1}{x^{\delta_2}} \exp[\frac{a_1 c_1}{1-\delta_1}x^{1-\delta_1}]$ is a monotonically   increasing   function.
       Then we obtain
     \begin{equation}
     \begin{split}
     & \sum_{s=0} ^k  \frac{1}{(s+1)^{\delta_2}} \exp[\frac{a_1 c_1}{1-\delta_1}(s+2)^{1-\delta_1}]\\
    &\leq \sum_{s=0} ^{k_1-1}  \frac{1}{(s+1)^{\delta_2}}  \exp[\frac{a_1 c_1}{1-\delta_1}(s+2)^{1-\delta_1}]
\\  &     \quad  +  \sum_{s=k_1} ^k  \frac{1}{(s/2+1)^{\delta_2}} \exp[\frac{a_1 c_1}{1-\delta_1}(s+2)^{1-\delta_1}]\\
  & =  \tau +  2^{\delta_2} \sum_{s=k_1+2} ^{k+2}  \frac{1}{s^{\delta_2}} \exp[\frac{a_1 c_1}{1-\delta_1}s^{1-\delta_1}]
       \\  &  \leq  \tau +  2^{\delta_2} \int_{x=k_1+2}^{k+3}
    \frac{1}{x^{\delta_2}} \exp[\frac{a_1 c_1}{1-\delta_1}x^{1-\delta_1}],
      \end{split}
     \end{equation}
     with  $\tau= \sum_{s=0} ^{k_1-1}  \frac{1}{(s+1)^{\delta_2}}  \exp[\frac{a_1 c_1}{1-\delta_1}(s+2)^{1-\delta_1}]$.
     It follows from  equation  \eqref{stepcomp} that
     \begin{equation}\label{stepcomp2}
     \begin{split}
      & \sum_{s=0} ^k \beta_{s} \mathbb{E} [\parallel \Phi(k,s+1 )\parallel ]   \leq
        a_2c_0 ( 2^{\delta_2} \frac{p(k)}{ q(k)} + \frac{\tau }{ q(k)}) ,
      \end{split}
     \end{equation}
     with $p(k)= \int_{x=k_1+2}^{k+3}  \frac{1}{x^{\delta_2}} \exp[\frac{a_1 c_1}{1-\delta_1}x^{1-\delta_1}]$.

  Even though both  $p(k)$ and $ q(k)$ tend  to infinity as $k \rightarrow \infty$,
 they are not increasing at the same order.
     In fact, we can the derivative of them as
     $$p(x)'=\frac{1}{(x+3)^{\delta_2}} \exp[\frac{a_1 c_1}{1-\delta_1}(x+3)^{1-\delta_1}],$$
     $$q(x)'= \frac{a_1c_1}{(x+2)^{\delta_1}} \exp[\frac{a_1 c_1}{1-\delta_1}(x+2)^{1-\delta_1}].$$

According to $\delta_2 >\delta_1$, $\lim_{x \rightarrow \infty } \frac{p(x)'}{q(x)'}= 0$, and hence, by
   the well-known L'H\^{o}pital's rule,  we get $\lim_{x \rightarrow \infty } \frac{p(x) }{q(x) }= 0$.
   Therefore, it follows from   \eqref{stepcomp2} that
    $$\sum_{s=0} ^k \beta_k \mathbb{E} [\parallel \Phi(k,s+1 )\parallel ] \rightarrow 0, k\rightarrow \infty. $$

(iii)
Note that $x^{-\delta_1}\exp[-c_3x^{1-\delta_2}]$ is a monotonically decreasing function when $c_3$ is positive,
and then with equations (\ref{fin}) and (\ref{stepestimate}), we have
     \begin{equation}\label{res2}
     \begin{split}
     & \sum_{k=0}^{\infty} \beta_{k+1} \mathbb{E} [\parallel \Phi(k,0) \parallel]  \\
     & \leq \sum_{k=0}^{\infty}  \frac{a_2}{(k+2)^{\delta_2}}
        c_0 \exp[\frac{a_1 c_1}{1-\delta_1}  -\frac{a_1 c_1}{1-\delta_1}(k+2)^{1-\delta_1}] \\
     & \leq c_7 \sum_{k=1}^{\infty}  \frac{1}{(k+1)^{\delta_1}} \exp[-c_8 (k+1) ^{1-\delta_1}]  \\
     & \leq  c_7 \int_{x=1}^{\infty}  \frac{1}{x^{\delta_1}} \exp[-c_8 x^{1-\delta_1}]\\
     & = - \frac{c_2}{c_8(1-\delta_1)}\exp[-c_8 x^{1-\delta_1}] \mid_{x=1}^{\infty} \\
     & =\frac{c_7}{c_8}\exp[-c_3  ],
      \end{split}
     \end{equation}
     with $c_7= a_2 c_0 \exp[\frac{a_1 c_1}{1-\delta_1}]$, $c_8= \frac{a_1 c_1}{1-\delta_1}$ as positive constants.

   (iv)
Since  $(x+2)^{- \eta} q(x)$ tends to infinity as $x \rightarrow \infty$,  with L'H\^{o}pital's rule,
      \begin{equation}
      \begin{split}
       &  \lim_{x \rightarrow \infty } \frac{p(x)}{ (x+2)^{-\eta}q(x)}\\
       &= \lim_{x \rightarrow \infty } \frac{p(x)'}{ (x+2)^{-\eta} q(x)' -\eta  (x+2)^{-\eta-1} q(x)} \\
       &= \lim_{x \rightarrow \infty } \frac{(x+2)^{-\delta_2}}{ a_1c_1(x+2)^{-\eta-\delta_1}-a_1c_1\eta  (x+2)^{-1-\eta}}\\
       &= \frac{1}{a_1c_1}. \nonumber
      \end{split}
     \end{equation}
Hence,  $$\frac{p(x)}{q(x)}= O( (x+2)^{-\eta})$$ with $\eta=\delta_2-\delta_1$.

     Therefore, by (\ref{stepcomp2}),
       \begin{equation}
       \begin{split}
       &\sum_{s=0} ^k \beta_s \mathbb{E} [\parallel \Phi(k,s+1 )\parallel]     \\
       & \qquad \leq  a_2c_0 ( 2^{\delta_2}  O( (k+2)^{-\eta})+ \frac{\tau}{q(k)}). \nonumber
       \end{split}
     \end{equation}
From $\delta_2 > \delta_1  $, we have $ \frac{1}{k^{\delta_2}} <  \frac{1}{k^{\delta_1}}$,
 and hence,   by a similar argument in \eqref{res2}  we get
$$\sum_{k=0}^{\infty}  \frac{1}{(k+2)^{\delta_2}}\exp[- c_3(k+2)^{1-\delta_1}]<\infty.$$
Finally,  we conclude that
     \begin{equation}\label{res}
     \begin{split}
       &\sum_{k=0}^{\infty} \beta_{k+1} \sum_{s=0} ^k \beta_{s} \mathbb{E} [\parallel \Phi(k,s+1 )\parallel]\\
       &\leq
        a_2^2c_0  2^{\delta_2}  \sum_{k=0}^{\infty}  \frac{O( (k+2)^{-\eta})}{(k+2)^{\delta_2}} + a_2^2c_0 \tau  \sum_{k=0}^{\infty} \frac{1}{(k+2)^{\delta_2}q(k)} \\
       & \leq  c_9 \sum_{k=0}^{\infty}  \frac{1}{(k+2)^{\delta_2}}\exp[- c_3(k+2)^{1-\delta_1}]\\
        & \qquad +c_{10} \sum_{k=0}^{\infty}  \frac{1}{(k+2)^{2 \delta_2-\delta_1}}   < \infty,
     \nonumber
     \end{split}
     \end{equation}
     with $c_9, c_{10}$ as  positive constants.
     \hfill $\Box$

\bibliographystyle{IEEEtran}
\bibliography{markov}

\end{document}